\documentclass{amsproc}

\usepackage{amssymb}

\usepackage{url}

\usepackage[dvips]{hyperref}

\let\mc=\mathcal \let\I=\mathbb
\let\ox=\otimes
  {\begin{list}{}%
      {\setlength\rightmargin{0pt}}%
    \item[]\it}
  {\end{list}}

\begin{document}

\title[Curve Counting]{Curve Counting  \`a la G\"ottsche}


\author[S. L. Kleiman]{Problem Session, 25 Aug 2011 --- updated, 13 Sep 2012\\
Steven L. Kleiman,\enspace MIT}
 \email{Kleiman@math.MIT.edu}
\thanks{The author offers his heartfelt thanks to
Florian Block, to Eduardo Esteves, to Ragni Piene, to Nikolay Qviller,
and to Vivek Shende for reading previous drafts and making apt comments.}

\subjclass[2010]{14N10 (Primary); 14C20, 14H20 (Secondary)}
\date{13 Sep 2012}

\maketitle

Let $n_\delta$ be the number of $\delta$-nodal curves lying in a
suitably ample complete linear system  $|L|$ and passing through appropriately
many points on a smooth projective complex algebraic surface.  Often
$n_\delta$ is referred to as a {\it Severi degree}.  A major
problem is to understand the behavior of $n_\delta$, specifically to
finish off Lothar G\"ottsche's mostly proved 1997 conjectures \cite{G97}
and then go on to treat the new refinements by G\"ottsche and Vivek
Shende \cite{G11}, \cite{GS12}.

The general area has been very active for over fifteen years, and is now
busier and more exciting than ever before.  Among many other people
involved have been Joe Harris himself, some of his students, and some of
theirs.  The area is unusually broad\,---\,embracing ideas from physics,
symplectic differential geometry, complex analytic geometry, algebraic
geometry, tropical geometry, and combinatorics.

{\bf Problem number one} is to find the two power series $$B_1(q),\,B_2(q)\in
\I Z[[q]]$$ appearing in G\"ottsche's remarkable formula for the
generating function of the $n_\delta$.  The formula expresses the
function, so the $n_\delta$, in terms of the four basic numerical
invariants of the system and the surface.  In fact, $n_\delta$ is a
polynomial in the four.  See (\ref{I}) and (\ref{P}) and (\ref{B})
below.

G\"ottsche \cite[Rmk.\,2.5(2)]{G97} computed the coefficients of
$B_1(q)$ and $B_2(q)$ up to degree $28$ on the basis of the recursive
formula for the $n_\delta$ of the plane due to Lucia Caporaso and Harris
\cite[Thm.\,1.1]{CH98}.  (A different recursion had been given earlier
by Ziv Ran \cite[Thm.\,3C.1]{R89}.)  G\"ottsche checked the result
against much of what was known, including Ravi Vakil's enumeration
\cite{V00} for the Hirzebruch surfaces (that is, the rational ruled
surfaces).

The  {\bf problem} is to find a closed
form for each $B_i(q)$, or a functional equation.

{\bf Second}, given $\delta$, how ample is suitable, so that $n_\delta$ has
the predicted value?  After all, for any system, the polynomial yields a
number, but it isn't always $n_\delta$.  For example, consider
plane curves of degree $d$.  If $d=1$, then $n_3$ is the number of
3-nodal lines, namely 0, but the polynomial yields 75.  Considering the
geometry, G\"ottsche \cite[Cnj.\,4.1, Rmk.\,4.4]{G97} conjectured the
polynomial works if  $d\ge\delta/2+1$.

The latter conjecture was proved for $\delta\le8$ by Ragni Piene and the
author \cite[Thm.\,3.1]{KP04} using algebraic methods, then for
$\delta\le14$ by Florian Block \cite[Prp.\,1.4]{B11}.  He built on ideas
of Sergey Fomin and Grigory Mikhalkin \cite[Thm.\,5.1]{FM10}, who used
tropical methods to set up the enumeration from scratch and to validate
its predictions for $d\ge2\delta$.  In principle, the problem is purely
combinatorial: to show formally the Caporaso--Harris recursion yields a
polynomial in $d$ for $d\ge\delta/2+1$.

On any surface, Martijn Kool, Shende, and Richard Thomas
\cite[Prp.\,2.1]{KST10} proved it suffices for $L$ to be $\delta$-very
ample.  Piene and the author \cite[Thm.\,1.1]{KP99} proved it suffices
for $L$ to be of the form $M^{\ox m} \ox N$ where $M$ is very ample,
$m\ge3\delta$, and $N$ is spanned, provided $\delta\le8$.  Both results
were inspired by G\"ottsche's \cite[Prp.\,5.2]{G97}; in turn, G\"ottsche
had been inspired by Harris and Rahul Pandharipande's paper \cite{HP95},
which treats the case $\delta\le3$ in the plane.

In fact, G\"ottsche conjectured the polynomial works for plane curves of
degree $d$ {\it iff\/} $d\ge\delta/2+1$.  And Block proved, for
$3\le\delta\le14$, that $\lceil\delta/2\rceil+1$ is, indeed, a {\it
threshold}, as he called it; namely, it is the least integer $d^*$ such
that the polynomial works for $d\ge d^*$.  Further, G\"ottsche
conjectured a similar statement for the Hirzebruch surfaces.  Shende and
the author \cite{KS12} proved that, above G\"ottsche's conjectured
threshold, the polynomials work for the plane and for the Hirzebruch
surfaces and that a similar statement holds for the classical del Pezzo
surfaces; moreover, there's at least one case where the polynomial 
works below the conjectured threshold too.

Sometimes, the curves are required to belong to a general linear
subsystem of $|L|$ rather than to pass through appropriately many
points.  However, the latter condition does yield a general subsystem by
Piene and the author's \cite[Lem.\,(4.7)]{KP99}.

The {\bf problem} is to determine just when the polynomial yields
$n_\delta$.

{\bf Third}, what about nonlinear systems?  After all, Gromov--Witten
theory fixes not the linear equivalence class, but the homology class,
and this class determines the four basic invariants, (\ref{I}) below.
Jim Bryan and Naichung Conan Leung \cite[Thm.\,1.1]{BL99} handled
primitive complete nonlinear systems on generic Abelian surfaces for all
$\delta$.  They used symplectic methods.  Piene and the author
\cite[\S\,5]{KP04} obtained similar results algebraically, but for
$\delta\le8$.

Israel Vainsencher \cite[\S\,6.2]{V95} treated a remarkable system.
His parameter space was the Grassmannian of $\I P^2$ in $\I P^4$.  His
surface was $\I P^2$, but moving in $\I P^4$.  His curves arose by
intersecting the moving $\I P^2$ with a fixed general quintic $3$-fold
$X$.  Thus he found $X$ contains 17,601,000 irreducible $6$-nodal
quintic plane curves.  Piene and the author \cite[Thm.\,4.3]{KP04}
validated the number.  Pandharipande \cite[(7.54)]{CK99} noted each
curve has six double covers previously unconsidered in mirror symmetry.

Given any suitably general algebraic system of curves on surfaces, 
Piene and the author \cite[Thm.\,2.5 and Rmk.\,2.7]{KP04} found  on
the parameter space the class of the curves with $\delta$ nodes for
$\delta\le8$ and conjectured the formula generalizes to any $\delta$. 

The {\bf problem} is to generalize the formula for
$n_\delta$, in (\ref{P}), to algebraic systems.

 {\bf Fourth}, what about higher singularities?  This question is
related to the previous one, about algebraic systems.  For example,
given a system, consider those curves with a triple point and $\delta$
double points.  Their number can be viewed as the number of curves with
$\delta$ double points in the following system: take the subsystem of
curves with a triple point, and resolve the locus of triple points.
This example was treated for $0\le\delta\le3$ by Vainsencher and by
Piene and the author \cite[Thm.\,1.2]{KP99}.  A substantial amount of
work has been done; see Maxim Kazarian's paper \cite{Ka03}, Dmitry
Kerner's papers \cite{K06}, \cite{K10}, Jun Li and Yu-Jong Tzeng's paper
\cite{LT12}, J\o{}rgen Rennemo's paper \cite{Re12} and their references.

The {\bf problem} is to enumerate the curves of fixed global
equisingularity type lying in a given algebraic system\,---\,that is, to
find on the parameter space the class of these curves.

{\bf Fifth}, what about positive characteristic?  Sometimes an enumeration is
more tractable modulo a prime.  Thus G\"ottsche \cite[Thm.\,0.1]{G90}
found the Betti numbers of the Hilbert schemes of points on a smooth
surface.  (In \cite[pp.\,175--178]{FG05}, he and Barbara Fantechi discuss
other proofs and refinements of the result.)  This result, combined with
others, led to the celebrated formula of Shing-Tung Yau and Eric Zaslow
\cite[p.\,5]{YZ96} enumerating rational curves on a K3 surface.  They
developed ideas of Cumrun Vafa et al.: see
\cite[p.\,438]{V96} for a similar formula; see \cite[p.\,44]{VW94}
for the use of G\"ottsche's result;  see \cite[p.\,437]{BVS96} for
the use of varying Jacobians.  In turn,
 Arnaud Beauville \cite{B99} and Fantechi, G\"ottsche, and Duco van
 Straten \cite{FGS99} developed the ideas in
 \cite{YZ96} further, and Xi Chen \cite[Thm.\,1]{Ch02} proved the curves
 are nodal.

The Yau--Zaslow formula too inspired G\"ottsche to develop his
conjectures.  For K3 surfaces and Abelian surfaces, $B_1(q)$ and
$B_2(q)$ disappear, leaving explicit formulas in any geometric genus.
These formulas were proved for primitive classes on generic such
surfaces by Bryan and Leung; see \cite{BL99b} for a lovely survey.

The {\bf problem} is to determine just when G\"ottsche's conjectures
hold in positive characteristic.

To define the $B_i(q)$, denote the surface by $S$, and its canonical
bundle by $K$.  The four basic invariants are these numbers:
\begin{equation}\label{I}
x:=L^2,\quad y:=L\cdot K,\quad z:=K^2,\quad t:=c_2(S).
  \end{equation}
For $\delta\le6$, Vainsencher \cite[\S\,5]{V95} worked out formulas for
the $n_\delta$, getting enormous polynomials in $x,y,z,t$.
Afterwards, it was natural to conjecture this statement:
\begin{equation}\label{V}
 \text{\it The number } n_\delta\text{ \it is given by a universal
   polynomial    of degree $\delta$ in } \I Q[x,y,z,t].
  \end{equation}

For plane curves of degree $d$, we have $(x,y,z,t)=(d^2,-3d,\,9,\,3)$.  So
Philippe Di Francesco and Claude Itzykson \cite[p.\,85]{dFI95}
conjectured $n_\delta$ is given by a polynomial in $d$ of a certain
shape for $\binom{d-1}2\ge\delta$.  Youngook Choi \cite[p.\,12]{C96}
established their conjecture for $d\ge\delta$ on the basis of Ran's work
\cite{R89}.  G\"ottsche \cite[\S\,4]{G97} refined the conjecture.  Given
(\ref{V}) in the form of (\ref{P}) below, Nikolay Qviller
\cite[\S\,4]{Q10} established most of G\"ottsche's refinements
concerning the shape.

In full generality, (\ref{V}) was given a symplectic proof and an
algebraic proof by Ai-ko Liu \cite{L00}, \cite{L04}.  It was recently
given new proofs by Tzeng \cite[Thm.\,1.1]{T10}
and Kool, Shende, and Thomas \cite[Thm.\,4.1]{KST10}; the former is
purely algebraic, whereas the latter also relies on topology.  These new
proofs have caused quite a stir!

G\"ottsche \cite[Cnj.\,2.1]{G97} did conjecture (\ref{V}) in full
generality, but his elaboration is far more important.  First, he proved
(\ref{V}) is equivalent to this statement:
\begin{equation}\label{A}
\sum n_\delta u^\delta = A_1^xA_2^yA_3^zA_4^t 
  \text{\it \quad for some\enspace }A_i\in\I Q[[u]].
\end{equation}

The $A_i$ are the exponentials of their logarithms.  Hence
(\ref{A}) is equivalent to this:
\begin{equation}\label{P}
n_\delta=P_\delta(a_1,\dotsc,a_\delta)/\delta!
 \text{\quad\it where\enspace}\textstyle
 \sum_{\delta\ge0}P_\delta u^\delta/\delta!
 = \exp\bigl(\sum_{\kappa\ge1} a_\kappa u^\kappa/\kappa!\bigr)
\end{equation}
for some \textit{linear} forms $a_\kappa(x,y,z,t)$.  The polynomials
$P_\delta$ were studied extensively in 1934 by Eric Temple Bell
\cite{B34}.  Piene and the author \cite[p.\,210]{KP99} determined
$a_\kappa$ for $\kappa\le8$, and found its coefficients to be integers.
Recently, Qviller \cite[Thm.\,2.4]{Q10} (see \cite[\S\,6]{Q11} too)
proved the coefficients are always integers.

The $B_i(q)$ appear in the next formula,  the
\textit{G\"ottsche--Yau--Zaslow Formula}:
 \begin{equation}\label{B}
\sum n_\delta u(q)^\delta
  = B_1(q)^zB_2(q)^yB_3(q)^\chi B_4(q)^{-\nu/2} 
\end{equation}
where $u(q),\,B_3(q),\,B_4(q)\in\I Z[[q]]$ are explicit quasi-modular forms
and where
$$\chi:=\chi(L)=(x-y)/2 +\nu\quad\text{and}\quad
  \nu:=\chi(\mc O_S) =(z+t)/12.$$
G\"ottsche \cite[Cnj.\,2.4]{G97} conjectured (\ref{B}).  He
\cite[Rmks.\,2.5(1), 3.1]{G97} noted (\ref{B}) implies (\ref{V}) and
generalizes the Yau--Zaslow Formula.  Tzeng \cite[Thm.\,1.2]{T10}
derived (\ref{B}) from
 (\ref{A}) via Bryan and Leung's work on K3
surfaces \cite[Thm.\,1.1]{BL00} and via Piene and the author's
\cite[Lem.\,5.3]{KP04}; the latter provides enough suitably
ample primitive classes.

Finally, G\"ottsche and Shende were inspired by Kool, Shende, and
Thomas's work to conjecture, {\it  among many other statements,} 
refinements \cite[Cnj.~75]{GS12} of the Caporaso--Harris and Vakil
recursions.  Further, G\"ottsche and Shende \cite[Cnj.\,5,\,7]{GS12}
refine (\ref{B}) with this conjecture: there should be polynomials
$n_{\delta}(v)\in \I Z[v]$ and power series with polynomial coefficients
$u(v,q),\,B_i(v,q)\in \I Q[v,v^{-1}][[q]]$ such that
\begin{equation*}\textstyle \sum  n_\delta(v)\, u(v,q)^\delta
 = B_1(v,q)^zB_2(v,q)^yB_3(v,q)^\chi B_4(v,q)^{-\nu/2}  \end{equation*}
 and such that putting $v=1$ recovers (\ref{B}).  Again $u(v,q)$ and
$B_3(v,q)$ and $B_4(v,q)$ are known; however, it is an {\bf open
problem} to find the geometric meaning of $n_\delta(v)$.

If $S$ is a real toric variety, then $n_\delta(-1)$ is conjectured in
\cite[Cnj.~90]{GS12} to be the tropical Welschinger invariant\,---\,the
number of real $\delta$-nodal curves lying in a suitably ample real
complete linear system and passing through a subtropical set of
appropriately many real points, each curve counted with an appropriate
sign.  The notion of subtropical set was introduced and studied by
Mikhalkin in \cite{M05}.  This conjecture is also stated by Block and
G\"ottsche in a paper currently being written; further, there the
conjecture is proved for $\delta\le8$ using methods like those in \cite{B11} 

The {\bf refined problem number one} is to find $B_1(v,q)$ and
$B_2(v,q)$.

\bibliographystyle{amsplain}

\end{document}